\documentclass[12pt]{article}
\makeatletter
\newcommand{\singlespacing}{\let\CS=\@currsize\renewcommand{\baselinestretch}{1}\tiny\CS}
\usepackage{epsfig}
\usepackage[utf8]{inputenc}
\usepackage[T1]{fontenc}
\usepackage{amsmath}
\usepackage{amsthm}
\usepackage{mathptmx}
\usepackage{authblk}
\usepackage{adjustbox}
\usepackage{float}
\usepackage{tgpagella}
\usepackage[parfill]{parskip}
\usepackage[colon,sort&compress,round,authoryear]{natbib}
\newcommand{\be}{\begin{equation}}
\newcommand{\ee}{\end{equation}}
\newcommand{\beanno}{\begin{eqnarray*}}

\newcommand{\eeanno}{\end{eqnarray*}}
\newcommand{\bea}{\begin{eqnarray}}
\newcommand{\eea}{\end{eqnarray}}
\newcommand{\ba}{\begin{array}}
\newcommand{\ea}{\end{array}}

\newcommand{\bc}{\begin{center}}
\newcommand{\ec}{\end{center}}

\newcommand{\ds}{\displaystyle}

\renewcommand{\baselinestretch}{1.49}
\newtheorem{theorem}{Result}[section]

\allowdisplaybreaks
\begin{document}
\title{\sc A New Decision Theoretic Sampling Plan for Exponential Distribution under Type-I Censoring}
\author{Deepak Prajapati\footnote{Department of Mathematics and Statistics, Indian Institute of Technology Kanpur, India} \& 
Sharmistha Mitra\footnote{Department of Mathematics and Statistics, Indian Institute of Technology Kanpur, India} \& 
Debasis Kundu\footnote{Department of Mathematics and Statistics, Indian Institute of Technology Kanpur, India, Corresponding author, 
e-mail: kundu@iitk.ac.in.}}

\date{}
\maketitle
\begin{abstract}
In this paper a new decision theoretic sampling plan (DSP) is proposed for Type-I censored exponential distribution.  The proposed DSP is based on a new estimator of the expected lifetime of an exponential distribution which always exists, unlike the usual maximum likelihood estimator.  The DSP is a modification of the Bayesian variable sampling plan of \cite{Lam:1994}.  
An optimum DSP is derived in the sense that it minimizes the Bayes risk.  In terms of the Bayes risks, it performs better than 
Lam's sampling plan and its performance is as good as the Bayesian sampling plan of \cite{LLH:2002}, although implementation of the DSP is very simple.  Analytically it is more tractable than the Bayesian sampling plan of \cite{LLH:2002}, and it can be easily generalized for any other loss functions also.
A finite algorithm is provided to obtain the optimal plan and the corresponding minimum Bayes risk is calculated. Extensive numerical comparisons with the optimal Bayesian sampling plan proposed by \cite{LLH:2002}  
are made.  The results have been extended for three degree polynomial loss function and for Type-I hybrid censoring scheme.

\end{abstract}

%\noindent$^1$ Department of Mathematics and Statistics, Indian Institute of Technology Kanpur, India.  \\
%\noindent$^2$ Corresponding author, e-mail: kundu@iitk.ac.in

\section{\sc Introduction} 
In statistical quality control, acceptance sampling plans play a crucial role in deciding whether to accept or reject a batch of products. There are various approaches to determine these sampling plans. Both classical and decision theoretic approaches have been discussed in the literature. From the point of view of economy, the decision theoretic approach is more scientific and realistic because the sampling plan is determined by making an optimal decision on the basis of some economic consideration such as maximizing the return or minimizing the loss.  Extensive work has been done dealing with the designing of the sampling plans under different censoring schemes.  See for 
example \cite{Hald:1967}, \cite{FM:1974}, \cite{Thyregod:1975},  \cite{WK:1979}, \cite{Lam:1988, Lam:1990}, 
\cite{LC:1995}, \cite{HL:2002}, \cite{LHB:2008a, LHB:2008b}, \cite{LY:2013}, \cite{TCLY:2014}, \cite{CYL:2015} and the references cited therein. 

In practice, life testing experiments are usually censored in the sense that the testing procedure terminates before the actual lifetime of all selected items are observed.   \cite{Lam:1994} investigated this problem of formulating sampling plans for the exponential distribution with Type-I censoring. The one-sided decision function is based on the maximum likelihood estimator (MLE) of the expected lifetime, when
it exists. The sampling plan is a triplet $(n,\tau,\xi)$ where $n$ is the number of items inspected, $\tau$ is the fixed censoring time and $\xi$ is the minimum acceptable surviving time, i.e., the lot is accepted if and only if the maximum likelihood estimate of the expected lifetime is larger than or equal to $\xi$. 

In this approach, there are two shortcomings which should be mentioned. First of all, under a Type-I censoring, the MLE of the expected lifetime 
may not always exist.  Secondly, the loss function of \cite{Lam:1994} consists of the sampling cost and a component of the polynomial decision loss function only.  \cite{LLH:2002} first observed that \cite{Lam:1994} does not take the cost of the testing time $\tau$ into account in the loss function. In fact, they rightly observed that taking $\tau=\infty$, the sampling plan proposed by them gives better results than the sampling plan of \cite{Lam:1994}. It indicates that if we do not consider the cost of the testing time in the loss function then it leads to wrong results.
By considering a Bayes decision function which minimizes Bayes risk over optimal $(n,\tau)$,  \cite{LLH:2002} showed that the 
sampling plan of \cite{Lam:1994} is ``neither optimal, nor Bayes''.  From now on the sampling plan proposed by  \cite{LLH:2002}, will be referred as LSP.

It should be noted that although the LSP is a Bayes sampling plan and it performs better
than the sampling plan of \cite{Lam:1994} in terms of lower Bayes risks, the implementation of the LSP may not be
very simple particularly if the loss function is complicated.   \cite{LLH:2002} provided a simple algorithm in explicit form when the loss function is a two degree polynomial.  In this paper our goal is to develop an optimal sampling plan which is 
easy to implement in practice even when the loss function is not simple, and which performs equally well as the LSP in terms of the Bayes risks.  We call this new sampling plan as the decision theoretic sampling 
plan (DSP).  An optimal DSP is derived in this case and numerical results show that it performs as good as the Bayesian sampling plan.  A finite algorithm is provided to obtain the optimal DSP.

\enlargethispage{1 cm}

We further consider the case when the loss function may not be a quadratic.  It is observed in this case that the proposed DSP is 
very easy to implement, where as to implement the LSP one needs to solve a higher degree polynomial. It is clear that the implementation of the LSP becomes more difficult 
if the loss function is a higher degree polynomial or if it is not a polynomial loss function, although the DSP can be implemented
quite conveniently for higher degree polynomial or for a non-polynomial loss function.  We present some numerical results for 
three degree polynomial and for non- polynomial loss function and it is observed that in terms of Bayes risks the optimal DSP is as good as the Bayesian sampling plan LSP. Finally, we have extended our results for Type-I hybrid censoring scheme also.

The rest of the paper is organized as follows.  In Section 2 we provide the DSP.  All the necessary theoretical results are provided
in Section 3. In Section 4 we discuss how DSP is more tractable then LSP in case of higher degree polynomial loss function.  In Section 5 we extend the results for hybrid Type-I censoring. Numerical results for optimal DSP and comparison between the DSP and LSP are 
provided in Section 6.  
Finally we conclude the paper in Section 7.

\section{\sc Problem Formulation and the Proposed Decision Rule}

Suppose $n$ identical items are put on a life testing experiment.  Let their lifetimes be denoted by $(X_1, X_2,\ldots, X_n)$ and let 
$X_{(1)}\leq X_{(2)}\leq \ldots \leq X_{(n)}$ be the corresponding order statistics of the random sample of size $n$.  We stop the experiment
at the time point $\tau$.  Hence, our observed sample is $(Y_1 ,\ldots,Y_M)$=$(X_{(1)},\ldots,X_{(M)})$, $M = \max \{i : X_{(i)}\leq \tau \}$ 
is the random number of failures that occurs before time $\tau$.  

\enlargethispage{1.5 cm}
It is assumed that the lifetime of an experimental unit follows exponential distribution with the parameter $\lambda$ and it has the 
following probability density function (PDF) 
\begin{equation}
f(x; \lambda) = 
\begin{cases}
\lambda e^{-\lambda x} & \mbox{if} \ x > 0 \\
0                        &  \mbox{otherwise.}
\end{cases}
\end{equation}
We denote $\theta = 1/\lambda$, the expected lifetime of an experimental unit.  Based on the random sample 
$Y = (Y_1, \ldots, Y_M)$, we define the decision function as
\begin{equation}
\delta(\textbf{y}) =
\begin{cases}
d_{0}         & \mbox{if}  \  \widehat{\theta} \geq \xi, \\
d_{1}         & \mbox{if}   \  \widehat{\theta} < \xi , 
\end{cases}    \label{dec-func}
\end{equation}
where $\widehat{\theta}$ is a suitably defined estimator of $\theta$ (need not be a MLE), $\xi$ denotes the minimum acceptable 
surviving time, $d_0$ and $d_1$ are decisions of accepting and rejecting the batch, respectively.
A polynomial loss function including cost of time is define as
\begin{equation}
L(\delta(\textbf{x}),\lambda,n, \tau)=
\begin{cases}
nC_s + \tau C_{\tau} + g(\lambda) &\mbox{if} \ \delta(\textbf{x}) = d_{0}, \\
nC_s + \tau C_{\tau} + C_r         &\mbox{if} \ \delta(\textbf{x}) = d_{1},
\end{cases}   \label{loss-func}
\end{equation}
here $ C_s, C_{\tau}, C_r $ are positive constants defining $C_s$ = inspection cost per item, $C_{\tau}$ = cost due to per unit time, 
$C_r$ = cost due to rejection of the batch, and 
$g(\lambda) = a_0 + a_1 \lambda + \ldots + a_k \lambda^k$ is the loss due to accepting the batch provided the coefficients 
$a_0, a_1, \ldots, a_k $ are such that 
\begin{equation}
g(\lambda) = a_0 + a_1 \lambda + \ldots + a_k \lambda^k \geq 0 \ \ \ \ \ \ \ \ \forall \ \  \lambda >0.
\end{equation}
Further, it is assumed that $\lambda$ has a conjugate gamma prior distribution with the shape parameter $a > 0$ and the scale parameter 
$b > 0$ denoted by $G(a,b)$ with the following PDF 
\begin{equation}
\pi(\lambda;a,b) = \left \{ 
\begin{array}{lll}
\frac{b^a}{\Gamma(a)} \lambda^{a-1} e^{-\lambda b} & \hbox{if} & \lambda > 0  \cr
0 & \hbox{if} & \lambda \le 0, \ \ 
\end{array}
\right .
\label{prior}
\end{equation}
where $a$ and $b$ are the hyper parameters.

Our aim is to determine the optimal sampling plan namely $(n, \tau, \xi)$ by using the decision function (\ref{dec-func}) so that the Bayes risk is minimized under the given loss function (\ref{loss-func}) over all such sampling plans. 
Let $\widehat{\theta}_{M}$ denote the maximum likelihood estimator (MLE) of the expected lifetime $\theta$, i.e.  
\begin{equation*}
\widehat{\theta}_M =  \left \{
\begin{array}{lll}
 \frac{\sum_{i=1}^{M} X_{(i)} +(n-M){\tau}}{M}  & \mbox{if} & M \geq 1 \cr 
 \hbox{does not exist} & \mbox {if} &    M=0.
\end{array}
\right .
\end{equation*}
\cite{Lam:1994} used the decision function (\ref{dec-func}) with 
\be
\widehat{\theta} = \left \{
\begin{array}{lll}
\widehat{\theta}_M & \hbox{if} & M \ge 1 \cr
n \tau & \hbox{if} & M = 0.  
\end{array}
\right .
\ee
But it is observed by  \cite{LLH:2002} that if there exists a cost function $C_{\tau} > 0$ in the loss function, then 
the decision function of \cite{Lam:1994} does not produce an optimal result.  

\enlargethispage{1 cm}
In this paper we propose a new estimator of $\theta$ for a given $c > 0$, 
and for all $M \ge 0$, as follows; 
\begin{equation}
\widehat{\theta}_{M+c}= \frac{\sum_{i=1}^{M} X_{(i)} + (n-M)\tau}{M+c} \ \ \ \ \forall \ \ \ M\geq 0.
\end{equation}
It is a shrinkage estimator of $\theta$ and we can see that for $c=0$ it is the  maximum likelihood of estimator $\theta$. It is straightforward to see that
$$
Var(\widehat{\theta}_{M+c} )\leq Var(\widehat{\theta}_{M}),
$$
and for $M=0$ this estimator exits when $c > 0$ so we use $\widehat{\theta} = \widehat{\theta}_{M+c}$ in (\ref{dec-func}).    
Any sampling plan with a given $n, \tau, \xi$ and $c$ is denoted by the quartet $(n, \tau, \xi, c)$.  The optimal sampling plan $(n_0, \tau_0, \xi_0, c_0)$ is obtained by 
determining that decision function in (\ref{dec-func}) with $\widehat{\theta} = \widehat{\theta}_{M+c}$, which minimizes the Bayes 
risk $r(n, \tau, \xi ,c)$ under the given loss function in (\ref{loss-func}) over all sampling plans $(n, \tau, \xi ,c)$.

\section{\sc Bayes Risks and Optimal Decision Rule}
\label{Sub3.1}
In this section we obtain the Bayes risks of the sampling plan $(n, \tau, \xi, c)$, and also the optimal decision rule $(n_0, \tau_0, 
\xi_0, c_0)$, which minimizes the Bayes risks.  In order to derive the Bayes risk of the sampling plan $(n, \tau, \xi, c)$, we first need to know the distribution of $\widehat{\theta}_{M+c}$.  It is clear that the distribution of $\widehat{\theta}_{M+c}$ has a discrete part and and an 
absolute continuous part.  The distribution of $\widehat{\theta}_{M+c}$ can be obtained as follows.
\begin{align}
\label{mix}
P(\widehat{\theta}_{M+c}\leq x) &= P(M=0)P(\widehat{\theta}_{M+c}\leq x|M=0)+P(M \geq 1)P(\widehat{\theta}_{M+c}\leq x|M \geq 1) \nonumber\\
&= p S(x)+ (1-p)H(x),
\end{align}
where $p=P(M=0)= e^{-n\lambda \tau}$,
\begin{eqnarray*}
S(x) & = & P(\widehat{\theta}_{M+c}\leq x|M=0)=
\left \{ 
\begin{array}{lll}
1  &\mbox{if} \ \ \  x \geq \frac{n\tau}{c}, \\
0   &\mbox{} \ \ \  otherwise,
\end{array} 
\right .
\\
H(x) & = & P(\widehat{\theta}_{M+c}\leq x|M \geq 1)=
\left \{
\begin{array}{lll}
\int_{0}^{x}h(u)du  &\mbox{if} \ \ \  0< x < \frac{n\tau}{c}, \\
0   &\mbox{} \ \ \  otherwise.
\end{array}
\right .
\end{eqnarray*}
Here $h(u)$ is the PDF of $\widehat{\theta}_{M+c}$ when $M \geq 1$ and 
\begin{align*}
h(u)=\frac{1}{1-p}\sum_{m=1}^{n}\sum_{j=0}^{m}\binom{n}{m}\binom{m}{j}(-1)^{j} e^{ -\lambda(n-m+j)\tau} \pi 
\big(u-\tau_{j; m, c};m,(m+c)\lambda \big)
\end{align*}
provided  $ 0< u < \frac{n\tau}{c}$, zero otherwise. To compute $h(u)$ for given $M\geq 1$, we use the similar approach as in 
\cite{GK:1998} and \cite{CCBK:2003}.  Alternatively, one can use B-spline technique provided in \cite{GC:2017} and 
\cite{CB:2013} to obtain $h(u)$.
 
\noindent Now we would like to compute the Bayes risk of the proposed sampling plan. To compare our results with those of  
\cite{LLH:2002}, we have assumed that \linebreak 
$g(\lambda)$ = $a_0 +a_1\lambda +a_2\lambda^2$ in (\ref{loss-func}), where $a_0,a_1$ and $a_2$ are fixed positive constants. Thus the loss function is defined as \\
\begin{equation}
L(\delta(\textbf{x}),\lambda,n, \tau)=
\begin{cases}
nC_s +  \tau C_{\tau} + a_0 +a_1\lambda +a_2\lambda^2  &\mbox{if} \quad \delta(\textbf{x}) = d_{0}, \\
nC_s + \tau C_{\tau} + C_r         &\mbox{if} \quad \delta(\textbf{x}) = d_{1}.
\end{cases}   \label{new-loss}
\end{equation}
The Bayes risk with respect to the loss function (\ref{new-loss}) is given by  
\begin{align}
r(n, \tau, \xi, c) = nC_{s}+ \tau C_{\tau}+a_0 + a_1 \mu_1 & +a_2 \mu_2  + \sum_{l=0}^{2}C_{l}\frac{b^a}{\Gamma (a)}\Bigg[\frac{\Gamma{(a+l)}}{(b+n\tau)^{(a+l)}} \textit{I}_{(\frac{n\tau}{c} < \xi)} 
 \nonumber \\            
 & + \sum_{m=1}^n\sum_{j=0}^m(-1)^{j}\binom{n}{m}\binom{m}{j} \frac{\Gamma{(a+l)}}{(C_{j,m})^{a+l}} 
I_{S^*_{j,m,c}}(m,a+l)\Bigg],     \label{brf}
\end{align}
where $\ds \mu_{i}= E(\lambda^i)$ for $i=1,2$ and $C_l$ is defined as 
\begin{equation}
C_l =
\begin{cases}
C_r- a_l  &\mbox{if} \quad l=0,   \nonumber \\
-a_l      &\mbox{if} \quad l= 1,2. \nonumber
\end{cases}
\end{equation}
The constant $C_{j,m}$, the exact expression of $I_{S^*_{j,m,c}}(m,a+l)$ and the proof of (\ref{brf}) are provided in the Appendix.

In general, for the loss function given in (\ref{loss-func}) of degree $k>2$, the Bayes risk can be evaluated in a similar way. 
\\
 \textit{Algorithm for finding Optimal sampling plan:}
 \\
To find the optimal value of sampling plan $n$, $\tau$, $\xi$ and $c$ based on the Bayes risk, a simple algorithm is described to obtain an optimal sampling plan in the following steps:
\begin{enumerate}
    \item Fix $n$ and $\tau$ minimize $r(n, \tau, \xi, c)$ with respect to $\xi$ and $c$ using grid search method  and denote the minimum Bayes risk 
     by $r(n,\tau,\xi_0 (n,\tau),c_0 (n,\tau))$.
    \item For fixed $n$ minimize  $ r(n,\tau,\xi_0 (n,\tau),c_0 (n,\tau))$ with respect to $\tau$ using grid search method and denote the minimum Bayes risk by $r(n, \tau_0 (n),\xi_0 (n,\tau_0 (n)),c_0 (n,\tau_0 (n)))$.
    \item Choose sample size $n_{0}$ such that  $$\ds r(n_0, \tau_0 (n_0), \xi_0 (n_0,\tau_0 (n_0)),c_0 (n_0,\tau_0 (n_0))) \leq r(n, \tau_0 (n),\xi_0 (n,\tau_0 (n)),c_0 (n,\tau_0 (n))) \  \ \forall \ \ \ n \geq 0 $$
\end{enumerate}
We denote the minimum Bayes risk by $r(n_0, \tau_0, \xi_0,c_0)$ with the corresponding optimal sampling plan $(n_0, \tau_0, \xi_0,c_0)$.
The following result implies that the algorithm is finite , i.e., we can find an optimal sampling plan in a finite number of search steps.
\begin{theorem}
 Assume $\xi$ has an upper bound since it is a minimum acceptable surviving time i.e $0<\xi<\xi^{*}$ and $0 < c < c^*$. Let us denote $r(n, \tau, \xi_0(n,\tau), c_0(n,\tau))
= min_{\xi, c} r(n, \tau, \xi, c)$ for some fixed $n \ (\geq 1)$ and $\tau$.  Further, let $n_0$ and $\tau_0$ be the optimal sample 
size and censoring time . Then,
\begin{align*}
n_0 \leq \min \bigg\{ \frac{C_r}{C_s},\frac{a_0 + a_1\mu_1 + \ldots + a_k \mu_k}{C_s},\frac{r(n, \tau, \xi_0(n,\tau), c_0(n,\tau))}{C_s}\bigg\} \\
\tau_0 \leq \min \bigg\{ \frac{C_r}{C_{\tau}},\frac{a_0 + a_1\mu_1 + \ldots + a_k \mu_k}{C_{\tau}},\frac{r(n, \tau, \xi_0(n,\tau), c_0(n,\tau))}{C_{\tau}}\bigg\}.
\end{align*}
\end{theorem}

\noindent {\bf Proof:} See in the Appendix.

\enlargethispage{1 in}

\section{\sc Higher degree polynomial loss function and Non polynomial loss function}
The main aim of this section is to show that if we have a higher degree polynomial loss function and for non polynomial loss function the implementation of DSP is much easier than LSP.
\subsection{\sc Higher degree polynomial loss function}
In general when the cost due to accepting the batch $g(\lambda)$ is a $k^{th}$ degree polynomial then the loss function is of the form given in (3) and corresponding Bayes risk expression of DSP is given by 
\begin{align}
r(n,\tau,\xi,c) = nC_{s}+\tau C_{\tau}+a_0& +a_1 \mu_1 + a_2 \mu_2 + \ldots + a_k \mu_k + \sum_{l=0}^{k}C_{l}\frac{b^a}{\Gamma (a)}\Bigg[\frac{\Gamma{(a+l)}}{(b+n\tau)^{(a+l)}} \textit{I}_{(\frac{n\tau}{c} < \xi)} 
 \nonumber \\            
 & + \sum_{m=1}^n\sum_{j=0}^m(-1)^{j}\binom{n}{m}\binom{m}{j} \frac{\Gamma{(a+l)}}{(C_{j,m})^{a+l}} 
I_{S^*_{j,m,c}}(m,a+l)\Bigg]. 
\label{DSP_H}
\end{align}
where $\mu_{i}= E(\lambda^i)$ for $i=1,2, \ldots, k$, $C_l$, $C_{j,m}$ and $ S^*_{j,m,c} $ have defined earlier.  For example, for a cubic polynomial loss function i.e $k=3$, the Bayes risk of DSP is given by
\begin{align}
r(n,\tau,\xi,c) = nC_{s}+\tau C_{\tau}+a_0 &+a_1 \mu_1 +a_2 \mu_2 + a_3 \mu_3 + \sum_{l=0}^{3}C_{l}\frac{b^a}{\Gamma (a)}\Bigg[\frac{\Gamma{(a+l)}}{(b+n\tau)^{(a+l)}} \textit{I}_{(\frac{n\tau}{c} < \xi)} 
 \nonumber \\            
 & + \sum_{m=1}^n\sum_{j=0}^m(-1)^{j}\binom{n}{m}\binom{m}{j} \frac{\Gamma{(a+l)}}{(C_{j,m})^{a+l}} 
I_{S^*_{j,m,c}}(m,a+l)\Bigg]. 
\end{align}
So for any integer $k$ the Bayes risk expression is simple and straightforward to calculate and we can also see that the form of the decision function is same for any value of $k$.
\par
Now in case of LSP, the Bayes decision function is given by
\begin{equation*}
\delta_B(\textbf{x}) =
\begin{cases}
1, & \mbox{if} \ \ \phi_{\pi}\big(m,y(n,\tau,m)\big) \leq C_{r}\\ 
0, &  \ \  \ \  \ \ otherwise.
\end{cases}
\end{equation*}
where 
\begin{align*}
&y(n,\tau,m)=\sum_{i=1}^{m} x_{i} + (n-m)\tau  \\
&\phi_{\pi}\big(m,y(n,\tau,m)\big) =\int_{0}^{\infty}g(\lambda)\pi\big(\lambda|m,y(n,\tau,m)\big)d\lambda.
\end{align*}
Since the prior distribution of $\lambda$ is $G(a,b)$, it immediately follows that the posterior distribution of $\lambda$ is 
also gamma and  
\begin{equation*}
\pi\big(\lambda|m,y(n,\tau,m)\big) \sim G(m+a,y(n,\tau,m)+b).
\end{equation*}
Then for a cubic loss function i.e $k=3$, 
\begin{align*}
    \phi_{\pi}\big(m, y(n,\tau,m)\big) & = \int_{0}^{\infty}g(\lambda)\pi\big(\lambda|m,y(n,\tau,m)\big)d\lambda\\
    &= a_{0} + \frac{a_{1}(m+a)}{(y(n,\tau,m)+b)}+ \frac{a_{2}(m+a)(m+a+1)}{(y(n,\tau,m)+b)^2} + \\
    & \ \ \ \ \ \ \ \ \  \frac{a_{3}(m+a)(m+a+1)(m+a+2)}{(y(n,\tau,m)+b)^3}.
\end{align*}
So to find the closed form of decision function we need to obtain the following set; 
$$
A = \{x;\ x \ge 0, \phi_{\pi}\big(m, x \big) \leq C_{r}\}.
$$
Observe that to construct the set $A$, we need to obtain the set of $x \ge 0$, such that   
\begin{equation}
h_1(x) = a_{0} + \frac{a_{1}(m+a)}{(x+b)}+ \frac{a_{2}(m+a)(m+a+1)}{(x+b)^2} + \frac{a_{3}(m+a)(m+a+1)(m+a+2)}{(x+b)^3}\leq C_{r},   \label{hdp}
\end{equation}
which is equivalent to find $x \ge 0$, such that  
\begin{eqnarray}
h_2(x) & = & (C_{r}-a_{0}) \big(x+b\big)^3-a_{1}(m+a)\big(x+b\big)^2   \nonumber \\ 
&  & -a_{2}(m+a)(m+a+1)\big(x+b\big)-a_{3}(m+a)(m+a+1)(m+a+2)\geq 0.   \label{hdp-2}
\end{eqnarray}
It can be easily shown that if $D_n(m)$ is the only real root or $D_n(m)$ is the maximum real root of $h_2(x)$ = 0, then the 
LSP will take the following form.
\begin{equation}
\delta_B(\textbf{x}) =
\begin{cases}
1, & \mbox{if} \ \  y(n,\tau,m) \geq D_{n}(m) - b \\ 
0, &  \ \  \ \  \ \ otherwise.
\end{cases}
\label{hdp-df}
\end{equation}
Therefore, to find the LSP,  one needs to solve a cubic equation which cannot be obtained in explicit form.  The associated Bayes
risk of (\ref{hdp-df}) can be obtained as given below;
\begin{equation*}
    r(n,\tau,\delta_B)=r_{1}+r_{2}+r_{3}+r_{4}\\
\end{equation*}
where 
\begin{align*}
    &r_{1}=nC_s + \tau C_{\tau}+ a_{0}+a_{1}\mu_{1}+a_{2}\mu_{2}+a_{3}\mu_{3}\\
    &r_{2}=I(n\tau<D_{n}(0)-b)\bigg\{\frac{(C_r-a_{0})b^{a}}{(n\tau+b)^a}-\frac{a_{1}ab^{a}}{(n\tau+b)^{a+1}}-\frac{a_{2}a(a+1)b^{a}}{(n\tau+b)^{a+2}}-\frac{a_{3}a(a+1)(a+2)b^{a}}{(n\tau +b)^{a+3}}\bigg\}\\
    &r_{3}=\sum_{m\in B}\sum_{j=0}^{l^*}\binom{n}{m}\binom{m}{j}\frac{(-1)^{j}}{(m-1)!}\frac{b^a}{\Gamma (a)}\big((n-m)\tau+b+j\tau\big)^{-(a+3)}\\  & \hspace{1cm}     \times \big\{(C_{r}-a_{0})\Gamma{(m+a)}\big((n-m)\tau+b+j\tau\big)^{3}\beta_{y}(m,a) \\ & \hspace{2cm} -a_{1}\Gamma{(m+a+1)}\big((n-m)\tau+b+j\tau\big)^{2}\beta_{y}(m,a+1)\\
    & \hspace{2.5cm}-a_{2}\Gamma{(m+a+2)}\big((n-m)\tau+b+j\tau\big)\beta_{y}(m,a+2)\\ & \hspace{5.5cm} - a_{3}\Gamma{(m+a+3)}\beta_{y}(m,a+3)\big\}\\
    &r_{4}=\sum_{m\in C}\sum_{j=0}^{m}\binom{n}{m}\binom{m}{j}(-1)^{j}b^{a}\big((n-m)\tau+b+j\tau\big)^{-(a+3)}\\
    & \hspace{1cm} \times \big\{(C_{r}-a_{0})\big((n-m)\tau+b+j\tau\big)^{3}-a_{1}a\big((n-m)\tau+b+j\tau\big)^{2}\\ & \hspace{2cm}-a_{2}a(a+1)\big((n-m)\tau+b+j\tau\big)-a_{3}a(a+1)(a+2)
    \big\}
\end{align*}
and \\
$ B=\{m\in I_{n}| 0<D_{n}(m)-b-(n-m)\tau\leq m\tau\}$, where $I_{n}=\{1,2, \ldots, n\}$ \\
$ C=\{m\in I_{n}| D_{n}(m)-b-(n-m)\tau>m\tau\}$\\
$y=\frac{(D_{n}(m)-b-(n-m)\tau-j\tau)}{D_{n}(m)}$\\
$l^{*}=\big[\frac{D_{n}(m)-b-(n-m)\tau}{\tau}\big]$, where $[x]$ largest integer not exceeding $x$.
\par
The problem becomes more complicated when $ k> 4 $ because it is well known that there is no algebraic solution to polynomial equations of degree five or higher (see chapter 5, \cite{Herstein:1975}. In these cases finding the Bayes risk for LSP is not straightforward but in the case of DSP it is quite easy as given in (\ref{DSP_H}).

\subsection{\sc  Non polynomial loss function}
Now consider the non-polynomial loss function where we will show that  the construction of LSP is not easy, where as DSP can be obtained quite easily. We consider a non polynomial loss function 
\begin{equation}
L(\delta(\textbf{x}),\lambda,n,\tau)=
\begin{cases}
nC_s + \tau C_{\tau} + a_0+a_1\lambda+a_2\lambda^{5/2} &\mbox{if} \ \delta(\textbf{x}) = d_{0}, \\
nC_s + \tau C_{\tau} + C_r         &\mbox{if} \ \delta(\textbf{x}) = d_{1},
\end{cases}   \label{nploss-func}
\end{equation}
where $g(\lambda)= a_0+a_1\lambda+a_2\lambda^{5/2}$ which is an increasing function in $\lambda$.
For the above proposed non-polynomial loss function the Bayes risk for DSP under Type-I censoring is given by
\begin{align}
r(n,\tau,\xi,c) = nC_{s}+\tau C_{\tau}+a_0 + a_1 & \mu_1 +a_2 \frac{\Gamma{(a+\frac{5}{2}})}{\Gamma{(a)}b^{\frac{5}{2}}} + \sum_{l=0}^{2}C_{l}\frac{b^a}{\Gamma (a)}\Bigg[\frac{\Gamma{(a+p_l)}}{(b+n\tau)^{(a+p_l)}} \textit{I}_{(\frac{n\tau}{c} < \xi)} 
 \nonumber \\            
 & + \sum_{m=1}^n\sum_{j=0}^m(-1)^{j}\binom{n}{m}\binom{m}{j} \frac{\Gamma{(a+p_l)}}{(C_{j,m})^{a+p_l}} 
I_{S^*_{j,m,c}}(m,a+p_l)\Bigg],
\end{align} 
where,
\begin{equation*}
p_l=
\begin{cases}
0, & \mbox{if} \ \ \ \ l=0, \\
1, & \mbox{if} \ \ \ \ l=1,\\
\frac{5}{2}, & \mbox{if} \ \ \ \ l=2.
\end{cases}
\end{equation*}
Now in case of LSP the Bayes decision function is given by 
\begin{equation*}
\delta_B(\textbf{x}) =
\begin{cases}
1, & \mbox{if} \ \ \phi_{\pi}\big(m,y(n,\tau,m)\big) \leq C_{r},\\ 
0, &  \ \  \ \  \ \ otherwise.
\end{cases}
\end{equation*}
where 
\begin{align*}
&y(n,\tau,m)=\sum_{i=1}^{m} x_{i} + (n-m)\tau  \\
&\phi_{\pi}\big(m,y(n,\tau,m)\big) =\int_{0}^{\infty}g(\lambda)\pi\big(\lambda|m,y(n,\tau,m)\big)d\lambda.
\end{align*}
In non polynomial loss function case $g(\lambda)=a_0+a_1\lambda+a_2\lambda^{5/2}$. So for the non polynomial loss function 
\begin{align*}
    \phi_{\pi}\big(m, y(n,\tau,m)\big) & = \int_{0}^{\infty}g(\lambda)\pi\big(\lambda|m,y(n,\tau,m)\big)d\lambda\\
    &= a_{0} + \frac{a_{1}(m+a)}{(y(n,\tau,m)+b)}+ \frac{a_{2}\Gamma(m+a+\frac{5}{2})}{\Gamma(m+a)(y(n,\tau,m)+b)^{\frac{5}{2}}}.
\end{align*}
So to find a closed form of the decision function we need to obtain the following set; 
$$
A = \{x;\ x \ge 0, \phi_{\pi}\big(m, x \big) \leq C_{r}\}.
$$
Observe that to construct $A$, we need to obtain the set of $x \ge 0$, such that
\begin{equation}
h_1(x) = a_{0} + \frac{a_{1}(m+a)}{(x+b)}+ \frac{a_{2}\Gamma(m+a+\frac{5}{2})}{\Gamma(m+a)(x+b)^{\frac{5}{2}}} \leq C_{r},   \label{hdnp} \nonumber
\end{equation}
which is equivalent to find $x \ge 0$, such that  
\begin{eqnarray}
h_2(x) & = & (C_{r}-a_{0})\Gamma(m+a) \big(x+b\big)^{\frac{5}{2}}-a_{1}(m+a)\Gamma(m+a)\big(x+b\big)^{\frac{1}{2}}-a_{2}\Gamma(m+a+\frac{5}{2})\geq 0.   \label{hdnp-2} \nonumber 
\end{eqnarray}
It is obvious that finding the closed form solution of the non polynomial equation $h_2(x)=0$ is not possible. So in case of the non-polynomial loss function it is difficult to construct the LSP and the explicit expression of the Bayes risk. But for the proposed method DSP this difficulty does not arises as DSP does not depend on the form of the loss function.

\section{\sc Type-I Hybrid censoring}
\label{Sub3.2}
When the random sample is coming from Type-I hybrid censoring.  Let us define  $\tau^{*}= min\{X_{(r)},\tau\}$ and $M^{*}$ is number of failure before time $\tau^{*}$. Then $M^{*}$ takes value $\{0,1,\ldots,r\}$ and for $M^{*}=0$ the MLE does not exist. We use new estimator which is define as,
$$ 
\widehat{\theta}_{M^{*}+c}  = \frac{\sum_{i=1}^{M^*}X_{(i)}+(n-M^{*})\tau^{*}}{M^{*}+c} 
                           = \begin{cases}
                           \frac{\sum_{i=1}^{M}X_{(i)}+(n-M)\tau}{M+c} & \mbox{if} \ \ \ \ \ \ \ X_{(r)}>\tau, \\
                           \frac{\sum_{i=1}^{r}X_{(i)}+(n-r)X_{(r)}}{r+c} & \mbox{if} \ \ \ \ \ \ \ X_{(r)}\leq, \tau,
                           \end{cases}
$$ 
where $M$ is a number of failure before time $\tau$.  In this case also 
\begin{align}
\label{mixh}
P(\widehat{\theta}_{M^{*}+c}\leq x) &= P(M^{*}=0)P(\widehat{\theta}_{M^{*}+c}\leq x|M^{*}=0)+P(M^{*} \geq 1)P(\widehat{\theta}_{M^{*}+c}\leq x|M^{*} \geq 1) \nonumber\\
&= p S(x)+ (1-p)H(x),
\end{align}
where $p=P(M^{*}=0)= e^{-n\lambda \tau}$, $S(x)$ and $H(x)$ are same as defined in (\ref{mix}) with the PDF $h(u)$ of $\widehat{\theta}_{M^{*}+c}$ when $M^{*}\geq 1$ is given as
\begin{align*}
h(u)= \frac{1}{1-p}\Bigg[\sum_{m=1}^{n}&\sum_{j=0}^{m}\binom{n}{m}\binom{m}{j}(-1)^{j} e^{ -\lambda \tau (n-m+j) }\pi 
\big(u-\tau_{j,m,c};m,(m+c)\lambda \big)+ \pi \Big(u ;r, r\lambda \Big)  \\
&
+r\binom{n}{r}\sum_{j=1}^{r}\binom{r-1}{j-1}\frac{(-1)^j e^{ -\lambda \tau (n-r+j) }}{n-r+j}\pi 
\big(u-\tau_{j,r,c};r,(r+c)\lambda \big) \Bigg]
\end{align*}
provided $0<u<\frac{n\tau}{c}$, zero otherwise. For computation of $h(u)$ for given $M\geq 1$, we have followed the method of 
\cite{GK:1998} and \cite{CCBK:2003} or we can use the approach proposed by \cite{GC:2017} and 
\cite{CB:2013}.  

\noindent Many recent research works on finding the Bayesian sampling plan is based on a quadratic  loss function (for example see, 
\cite{Lam:1990, Lam:1994}, \cite{LC:1995}, \cite{LLH:2002}, \cite{HL:2002},  \cite{LHB:2008a, LHB:2008b, LHB:2011}, \cite{LY:2013} etc). They used this functional form because computation become easier and $g(\lambda)$ = $a_0 +a_1\lambda +a_2\lambda^2$ is an approximation of true acceptance cost. However, it is well known that a higher degree polynomial is a better approximation of the true acceptance cost. So for better approximation we consider the functional form of the loss function defined as
\begin{equation}
L(\delta(\textbf{x}),\lambda)=
\begin{cases}
nC_s -(n-M)r_{s}+ \tau^{*} C_{\tau} + a_0 +a_1\lambda +\ldots+a_k\lambda^k  &\mbox{if} \ \delta(\textbf{x}) = d_{0}, \\
nC_s -(n-M)r_{s}+ \tau^{*} C_{\tau} + C_r         &\mbox{if} \ \delta(\textbf{x}) = d_{1},
\end{cases}   \label{loss-func1}
\end{equation}
where decision function $\delta(\textbf{x})$ is define in (\ref{dec-func}) with $\widehat{\theta}= \widehat{\theta}_{M^{*}+c}$. Using the 
distribution function (\ref{mixh}) and the loss function (\ref{loss-func1}), the Bayes risk of the DSP $(n, r, \tau, \xi,c)$ is computed similarly as in Section 3 and by \cite{LHB:2008a}.
\begin{theorem}
\label{T_3.3}
The Bayes risk of DSP $(n, r, \tau, \xi,c)$ w.r.t loss function (\ref{loss-func1}) is given as follows

\bea
 r(n, r,\tau, \xi, c) & = & n(C_{s}-r_{s}) + E(M^{*})r_{s} + E(\tau^{*}) C_{\tau}+a_0 +a_1 \mu_1 +\ldots+a_k \mu_k \nonumber\\
            &  & + \sum_{l=0}^{k}C_{l}\frac{b^a}{\Gamma (a)}\bigg\{\frac{\Gamma{(a+l)}}{(b+n\tau)^{(a+l)}} \textit{I}_{(\frac{n\tau}{c} < \xi)}+
 \sum_{m=1}^{n}\sum_{j=0}^{m}\binom{n}{m}\binom{m}{j}(-1)^j R_{l,j,m} \nonumber \\
            &  & \hspace{2.5cm}   + R_{l,r-n,r} + \sum_{j=1}^{r}\binom{n}{r}\binom{r-1}{j-1}(-1)^j \frac{r}{(n-r+j)} R_{l,j,r}\bigg\}, \nonumber
\eea
\noindent where
\bea
R_{l,j,m} & = &  \frac{\Gamma{(a+l)}}{(C_{j,m})^{a+l}} 
I_{S^*_{j,m,c}}(m,a+l), \label{rjm} \nonumber \\
E(M^{*})& = &\sum_{m=1}^{r-1}\sum_{j=0}^{m}m\binom{n}{m}\binom{m}{j}(-1)^{j}\frac{b^a}{(b+(n-m+j)\tau)^a}, \nonumber \\ 
& & \hspace{2.5cm}+\sum_{k=r}^{n}\sum_{j=0}^{k}r\binom{n}{k}\binom{k}{j}(-1)^j\frac{b^a}{(b+(n-k+j)\tau)^a}, \nonumber \\
E(\tau^{*})& = & r\binom{n}{r}\sum_{j=0}^{r-1}\binom{r-1}{j}(-1)^{r-1-j}\bigg\{\frac{b}{(n-j)^2(a-1)}- \frac{tb^a}{(n-j)((n-j)\tau+b)^a}, 
 \nonumber \\
&  & -\frac{b^a}{(n-j)^2(a-1)((n-j)\tau+b)^{a-1}}\bigg\} +  \sum_{k=r}^{n}\sum_{j=0}^{k}\tau\binom{n}{k}\binom{k}{j}(-1)^j\frac{b^a}{(b+(n-k+j)\tau)^a}. \nonumber 
\eea
For computation of $E(M^{*})$ and $E(\tau^{*})$ see \cite{LY:2013}.
\end{theorem}

\noindent Based on the explicit expression of the Bayes risk, an optimum DSP $(n_0, r_0, \tau_0, \xi_0,c_0)$ can be determined by 
\bea
r(n_0, r_0, \tau_0, \xi_0, c_0)=\underset{n,r\leq n}{min}\{\underset{\tau}{min}\{\underset{\xi,c}{min}[r(n, r, \tau, \xi, c)]\}\}.
\eea
\noindent In this case also the Bayes risk expression is very complicated so a similar algorithm as given in Section \ref{Sub3.1} is consider to obtain the optimum DSP $(n_0, r_0, \tau_0, \xi_0,c_0)$.
\par 
For each sample size $n$ and for given value of $r$, $\xi$ and $c$, the  Bayes risk  $r(n, r,\tau, \xi,c)$ is a function of $\tau$. If we have to find the minimum Bayes risk, we need an upper bound of $\tau$. \cite{TCLY:2014} suggested to choose a suitable range of $\tau$, say $[0, \tau_\alpha]$ where $\tau_\alpha$ is such that $P(0<X<\tau_{\alpha})=1-\alpha$, and $\alpha$ is preassigned number satisfying $0<\alpha<1$. The choice of $\alpha$ depend on the prescribed precision. The higher the precision required, the smaller the value of $\alpha$ 
should be taken.  Here we have used $\alpha=0.01$. In the range $[0, \tau_\alpha]$,  we have used grid search method to find the 
optimal value of $\tau$. Next result shows that the algorithm is finite , i.e., we can find an optimal sampling plan in a finite number of search steps.  The proof of the following result can be obtained similarly as the proof of {\bf Result 3.1}.
\newpage
\begin{theorem}
\label{T_3.4}
 Assuming that  $0<\xi \leq \xi^{*}$ and $0<c \leq c^{*}$ . Let us denote $r(n, r,\tau, \xi',c')$
= $\underset{\xi,c}{min} \  r(n,r,\tau, \xi,c)$ for some fixed $n \ (\geq 1)$ and $\tau$.  Further, let $n_0$ be the optimal sample 
size. Then,
\begin{align*}
&n_0 \leq \min \bigg\{ \frac{C_r}{C_s-r_s},\frac{a_0 + a_1\mu_1 + \ldots + a_k \mu_k}{C_s-r_s},\frac{r(n, r,\tau, \xi',c')}{C_s-r_s}\bigg\}
\end{align*}
and $r_{0}\leq n_{0}$.
\end{theorem}

\noindent The number of grid search points we choose depends on how well the behavior of the Bayes risk function $r(n,\tau, \xi,c)$ or $r(n, r, \tau, \xi,c)$ is. In practice,  if the values of Bayes risk function $r(n, \tau, \xi,c)$ or $r(n, r, \tau, \xi,c)$ is monotone or has unique minimum (in numerical examples we show this property) we will use less numbers of grid points. If the Bayes risk function  $r(n, \tau, \xi,c)$ or $r(n, r, \tau, \xi,c)$ are such as two or more sampling plans give the values equal to or close to the minimum value, then more grids search point are used and the grid search algorithm needs to be modified appropriately.

\section{\sc Numerical Results and Comparisons} 
 For the quadratic loss function to obtain the numerical results we consider the algorithm proposed in Section 3. Since the expression of  $r(n, \tau, \xi,c)$ is quite complicated, so it is not possible to obtain the optimal value of $n$, $\tau$, $\xi$ and  $c$ analytically. We need to obtain the optimal values of $n$, $\tau$, $\xi$ and  $c$ numerically so we need a following algorithm for that purpose:
 \newline
 \textit{Step-1:}
 For fixed $n$ and $\tau$, find the optimal values of $\xi$ and $c$, $\xi_0(n,\tau)$, $c_0(n,\tau)$, respectively, using grid 
search method.  The grid sizes of $\xi$ and $c$ are $0.0125$ and $0.0025$, respectively.
 \newline
 \textit{Step-2:}
 Let $n^{*}=min (C_r ,a_0 + a_1\mu_1+ a_2 \mu_2, r(n, \tau, \xi_0(n,\tau), c_0(n,\tau))) /\ C_s $ and $\tau_{*}=min(C_r,a_0 + a_1\mu_1+ a_2 \mu_2,r(n, \tau, \xi_0(n,\tau), c_0(n,\tau)))/\ C_{\tau}$ then it is clear that both $n^{*}$ and $\tau_{*}$ are finite and from the Result 3.1 ,
 $
 0\leq n_0 \leq n^{*} \ \ \ and \ \ \  0 \leq \tau_0 \leq \tau_{*} 
 $
Next for each $n$, compute  $r(n,\tau(n),\xi_0(n,\tau),c_0(n,\tau))$ and  minimize $r(n,\tau(n),\xi_0(n,\tau),c_0(n,\tau))$ with respect to $\tau$ where grid point are taken for $\tau$ is \ $0(0.0125)\tau_{*}$. Let the minimizer be denoted by $\tau_0(n)$.  
\newline
 \textit{Step-3:}
Finally choose that n for which  $r(n,\tau_0(n),\xi_0(n,\tau_0),c_0(n,\tau_0))$ is minimum.

\begin{table}
\caption{Minimum Bayes risk and corresponding optimal sampling plan for different values of $a$ and $b$ }
\label{table-1}
	\vspace{1em}	
	\begin{center}
		\begin{tabular}{|cc|ccccc|}\hline
	   $a$ & $b$ &
	   $r( n_0, \tau_0, \xi_0, c_0)$ & $n_0$ & $\tau_0$ & $\xi_0$ & $c_0$ \\ \hline
		0.2 & 0.2 & 9.0726 & 2 & 0.4625 & 0.2000 & 0.9600 \\ \hline
	    1.5 & 0.8 & 16.8439 & 3 & 0.4750 & 0.2250 & 0.1100 \\ \hline
        2.0 & 0.8 & 21.5046 & 3 & 0.6000 & 0.2750 & 0.1025 \\ \hline
		2.5 & 0.6 & 28.1949 & 3 & 0.8625 & 0.3125 & 0.8650 \\ \hline
		2.5 & 0.8 & 25.2777 & 3 & 0.7250 & 0.3000 & 0.3550 \\ \hline
		2.5 & 1.0 & 22.0361 & 3 & 0.5625 & 0.2625 & 0.0725 \\ \hline
		3.0 & 0.8 & 28.0087 & 3 & 0.8250 & 0.3125 & 0.7125 \\ \hline
		3.5 & 0.8 & 29.7131 & 2 & 0.8125 & 0.4125 & 0.4400 \\ \hline
		10.0 & 3.0 & 29.8053 & 1 & 0.4375 & 0.4750 & 0.8075 \\ \hline
		\end{tabular}
  \end{center}
\end{table}
We present some DSP for different values of $a$ and $b$ in Table 1. We have taken the following configuration:
$ a_0 = 2, a_1 = 2, a_2 = 2, C_r = 30, C_s = 0.5, C_{\tau} = 0.5, \xi^{*}=2 , c^{*}=1. 
$
Here $r(n_0, \tau_0, \xi_0,c_0)$ denotes the minimum Bayes risk, while $(n_0,\tau_0,\xi_0,c_0)$ is the corresponding optimal sampling plan. 
For example, consider the parameter values $a=2.5,\ b=0.8$ and coefficients $\ a_0=2,\ a_1=2,\ a_2=2,\ C_s=0.5,\ C_{\tau}=0.5,\ 
C_r=30$, the minimum Bayes risk is $r(3,0.7250,0.3000,0.3550)=25.2777$ indicating that the corresponding optimal sampling plan is 
$(3,0.7250,0.3000,0.3550)$ as given in Table \ref{table-1}. Thus, if we take 3 items from a batch to test under Type-I censoring, 
at a censoring time of 0.7250, we may accept that batch if the estimated average lifetime of the items $\hat{\theta}_{M+c}$  is greater than or equal to 0.3000 and the value of $c_0=0.3550$ ensures the existence of such an estimator.
\par
From Table 1 we can see that for a fixed value of $b$ as we increase the value of $a$, the Bayes risk increases. However, for a fixed value of $a$  as we increase the value of $b$, the Bayes risk decreases. At the same time, it is also seen that as the shape parameter of the prior distribution ($a$) increases, the minimum Bayes' risk increases irrespective of whether its scale parameter ($b$) increases or decreases.
\subsection{\sc Comparison between the LAM's sampling plan and the DSP}

In this section, we focus on comparison between the sampling plan of \cite{Lam:1994} and the DSP, some numerical results are presented in Table \ref{table-lam}. The values of coefficients $\ a_0=2,\ a_1=2,\ a_2=2,\ C_s=0.5,\ C_{\tau}=0$ and $C_r=30$ are used for comparison. In Table \ref{table-lam} only hyper-parameter $a$ and $b$ are varying and others are kept fixed.
\begin{table}[H]
\caption{Minimum Bayes risk and corresponding optimal sampling plan for different values of $a$ and $b$}
\label{table-lam}
	\vspace{1em}	
	\begin{center}
		\begin{tabular}{|c|cc|ccccc|}\hline
	$M$ &  $a$ & $b$ & $r( n_0, \tau_0, \xi_0, c_0)$ & $n_0$ & $\tau_0$ & $\xi_0$ & $c_0$ \\ \hline
	DSP  & 0.2 & 0.2 & 8.8228 & 2 & 0.6000 & 0.1875 & 1.1575 \\ 
	LAM  &     &     & 12.1499 & 4 & 0.0270 & 0.1080 &  \\ \hline
	DSP & 1.5 & 0.8 & 16.5825 & 3 & 0.7000 & 0.1750 & 1.0000 \\ 
	LAM  &     &     & 16.6233 & 3 & 0.5262 & 0.2631 &  \\ \hline
    DSP & 2.0 & 0.8 & 21.1398 & 4 & 1.1625 & 0.2000 & 1.7975 \\ 
    LAM  &     &     & 21.2153 & 3 & 0.6051 & 0.3026 &  \\ \hline
    DSP & 2.5 & 0.4 & 29.7506 & 1 & 0.8000 & 0.3250 & 1.4400 \\ 
    LAM  &     &     & 29.7506 & 1 & 0.7978 & 0.7978 &  \\ \hline
	DSP & 2.5 & 0.6 & 27.7266 & 3 & 1.2125 & 0.2750 & 1.3875 \\ 
	LAM  &     &     & 27.7834 & 3 & 0.8537 & 0.4268 &  \\ \hline
	DSP & 2.5 & 0.8 & 24.8419 & 4 & 1.3125 & 0.3000 & 0.3750 \\ 
	LAM  &     &     & 24.9367 & 3 & 0.7077 & 0.3539 &  \\ \hline
	DSP & 2.5 & 1.0 & 21.7081 & 4 & 1.1125 & 0.2250 & 0.9450 \\ 
	LAM  &     &     & 21.7640 & 3 & 0.5483 & 0.2742 &  \\ \hline
	DSP & 3.0 & 0.8 & 27.5581 & 3 & 1.1625 & 0.3000 & 0.8650 \\ 
	LAM  &     &     & 27.6136 & 3 & 0.8170 & 0.4085 &  \\ \hline
	DSP & 3.5 & 0.8 & 29.2789 & 2 & 1.0125 & 0.2750 & 1.6600 \\ 
	LAM  &     &     & 29.2789 & 2 & 1.0037 & 0.5019 &  \\ \hline
	DSP & 10.0 & 3.0 & 29.5166 & 2 & 0.8000 & 0.2625 & 1.0250 \\ 
	LAM  &     &     & 29.5166 & 2 & 0.7928 & 0.3964 &  \\ \hline
		\end{tabular}
  \end{center}
\end{table}
From Table-\ref{table-lam}  it is clear that the Bayes risk of the optimal DSP is less then or equal to the Bayes risk of Lam's sampling plan. Therefore, the optimal DSP is a better sampling plan then the Lam's sampling plan.

\subsection{\sc Comparison between the LSP and the DSP }
\subsubsection{\sc Comparison in terms of Bayes risk under quadratic loss function}
In this section, we present some numerical results to compare DSP and LSP. We have taken the same loss function  as in  
\cite{LLH:2002}  where coefficients $\ a_0=2,\ a_1=2,\ a_2=2,\ C_s=0.5,\ C_{\tau}=0.5 $ and $C_r=30$ are used for comparison. The results are presented in Table \ref{table-2}. 
 
\begin{table}[H]
\caption{Minimum Bayes risk and corresponding optimal sampling plan for different values of $a$ and $b$ and for DSP and LSP}  \label{table-2}
	\vspace{1em}
	\begin{center}
		\begin{tabular}{|cc|c|ccccc|}\hline
		\multicolumn{1}{|c}{$a$} &
		\multicolumn{1}{c|}{$b$} &
		\multicolumn{1}{c|}{LSP} &
		\multicolumn{5}{|c|}{DSP} \\ \cline{3-8}
	   & & $r( n_B, \tau_B, \delta_B)$ &
	   $r( n_0, \tau_0, \xi_0, c_0)$ & $n_0$ & $\tau_0$ & $\xi_0$ & $c_0$  \\ \hline
		0.1 & 0.2 & 6.1832 & 6.1832 & 2 & 0.4000 & 0.2000 & 0.8050 \\ \hline
		1.0 & 0.2 & 24.8966 & 24.8966 & 3 & 0.8250 & 0.3125 & 0.6700 \\ \hline
		1.5 & 0.8 & 16.8439 & 16.8439 & 3 & 0.4750 & 0.2250 & 0.1100  \\ \hline
		1.5 & 2.0 & 5.3750 & 5.3750  & 0 & 0      & 0      & 0        \\ \hline
		2.5 & 0.8 & 25.2777 & 25.2777 & 3 & 0.7250 & 0.3000 & 0.3550 \\ \hline
		2.5 & 1.0 & 22.0361 & 22.0361 & 3 & 0.5625 & 0.2625 & 0.0725  \\ \hline
		2.5 & 1.2 & 18.3194 & 18.3194 & 0 & 0      & 0      & 0        \\ \hline
		3.0 & 0.8 & 28.0087 & 28.0087 & 3 & 0.8250 & 0.3125 & 0.7125   \\ \hline
		3.5 & 0.8 & 29.7131  & 29.7131 & 2 & 0.8125 & 0.4125 & 0.4400   \\ \hline
		\end{tabular}
  \end{center}
\end{table}
\noindent where $r(n_0, \tau_0, \xi_0, c_0)$ denotes the minimum Bayes risk while $(n_0,\tau_0,\xi_0,c_0)$ is the corresponding optimal 
sampling plan for DSP, whereas $r(n_B, \tau_B, \delta_B)$ denote Bayes risk for LSP. From Table \ref{table-2} it is observed that in terms of Bayes risk of optimal DSP is good approximation of LSP.  It is also observed that for certain set of values of the hyper parameters and costs of the loss function, 
the optimal sample size and the censoring time are $0$ for both the plans.  It means the decision rule suggests acceptance of the
lot without any inspection in such cases.

\subsubsection{\sc Comparison in terms of Proportion of Acceptance under quadratic loss function}
For some further analysis we will give proportion of acceptance of some selected optimal sampling plans. We taken the following set of coefficients $\ a_0=2,\ a_1=2,\ a_2=2,\ C_s=0.5,\ C_{\tau}=0.5,\ C_r=30$ and parameter values $a=2.5,\ b=0.8$ . The results are presented in Table 3 by varying $a$ and $b$ and keeping other fixed . In Table 4, we have reported the results for different values of $a_0$, $a_1$ and $a_2$. Similarly, in Table 5, we have reported the results for different values of $C_s$ ,$C_\tau$ and $C_r$. All the results are based on $10000$ replications. In all the cases the proportion of acceptance are very high for both DSP and LSP and they are very close to each other.
\begin{table}[H]
\caption{Proportion of acceptance of  DSP and LSP for different values of $a$ and $b$.}  \label{table-3}
\vspace{1em}
	\begin{center}
		\begin{tabular}{|cc|c|c|}\hline
	   $a$ & $b$ & DSP & LSP \\ \hline
		1.7 & 0.2 & 0.8440 & 0.8424  \\ \hline
		2.1 & 0.3 & 0.7428 & 0.7443 \\ \hline
		2.4 & 0.4 & 0.7414 & 0.7262 \\ \hline
		\end{tabular}
  \end{center}
\end{table}
\begin{table}[H]
\caption{Proportion of acceptance of  DSP and LSP for different values of $a_0$, $a_1$ and  $a_2$.}  \label{table-4}
\vspace{1em}
	\begin{center}
		\begin{tabular}{|c|c|c|c|c|c|c|c|c|}\hline
	   $a_0$ & DSP & LSP & $a_1$ & DSP & LSP & $a_2$ & DSP & LSP \\ \hline
		13.5 & 0.7348 & 0.7219 & 10.2 & 0.6176 & 0.6122 &  6.0 & 0.8595 & 0.8624\\ \hline 
		14.0 & 0.7170 & 0.7096 & 10.5 & 0.6193 & 0.5956 & 6.5 & 0.7662 & 0.7686\\ \hline 
		14.5 & 0.8198 & 0.8134 & 10.8 & 0.5981 & 0.5822 &  6.8 & 0.7661 & 0.7586 \\ \hline
		\end{tabular}
  \end{center}
\end{table}
\begin{table}[H]
\caption{Proportion of acceptance of  DSP and LSP for different values of $C_s$, $C_{\tau}$ and $C_r$.}  \label{table-5}
\vspace{1em}
	\begin{center}
		\begin{tabular}{|c|c|c|c|c|c|c|c|c|}\hline
		 $C_s$ & DSP & LSP & $C_{\tau}$ & DSP & LSP & $C_r$ & DSP & LSP \\ \hline
		 3.0 & 0.9309 & 0.9354 & 3.0 & 0.9980 & 0.9976 &  17.5 & 0.7494 & 0.7389  \\ \hline  
		 3.5 & 0.8911 & 0.8988 & 3.5 & 0.9988 & 0.9988 &  18.0 & 0.7180 & 0.7081 \\ \hline
		 4.0 & 0.8913 & 0.8981 & 4.0 & 0.9981 & 0.9980 &  18.5 & 0.8030 & 0.8008 \\ \hline
		\end{tabular}
  \end{center}
\end{table}

\subsubsection{\sc Comparison in terms of Bayes risk under Higher degree polynomial and Non polynomial  loss function}
In Section 4 we have developed the theoretical results for higher degree polynomial and for non polynomial loss function. Where  we have shown that implementation of DSP is quite easy compare to LSP. Now to compare the performances of DSP and LSP for higher degree polynomial loss function we take cubic polynomial loss function  and we consider the following coefficients:
$$ 
a_0 = 2,a_1 = 2,a_2 = 2,a_3=2 , C_r = 30, C_s = 0.5, C_\tau = 0.5,  \xi^{*}=2, c^{*}=2.
$$
\enlargethispage{1.00in}
\begin{table}[H]
\caption{Minimum Bayes risk and corresponding optimal sampling plan for different $a$ and $b$ for DSP and LSP for cubic loss function}
	\vspace{1em}
	\begin{center}
		\begin{tabular}{|cc|c|ccccc|}\hline
		\multicolumn{1}{|c}{$a$} &
		\multicolumn{1}{c|}{$b$} &
		\multicolumn{1}{c|}{LSP} &
		\multicolumn{5}{|c|}{DSP}  \\ \cline{3-8}
	   & &
	   $r( n_B, \tau_B, \delta_B)$ &
	   $r( n_0, \tau_0, \xi_0, c_0)$ & $n_0$ & $\tau_0$ & $\xi_0$ & $c_0$  \\ \hline
		0.1 & 0.2 & 7.4606 & 7.4606 & 2 & 0.8875 & 0.3500 & 1.4875 \\ \hline
		0.5 & 0.8 & 10.0670 & 10.0670 & 3 & 0.8500 & 0.4250 & 0.0875  \\ \hline
		1.0 & 0.2 & 27.6919 & 27.6919 & 3 & 1.3625 & 0.5125 & 1.2750  \\ \hline
		1.0 & 0.8 & 17.0625 & 17.0625 & 4 & 1.1375 & 0.5000 & 0.1750 \\ \hline
		1.5 & 0.8 & 22.9149 & 22.9149 & 4 & 1.3000 & 0.5000 & 0.6875  \\ \hline
		2.5 & 0.8 & 29.7994 & 29.7994 & 2 & 1.4500 & 0.5750 & 1.2000 \\ \hline
		2.5 & 1.0 & 28.2333 & 28.2333 & 4 & 1.3250 & 0.5000 & 1.2875  \\ \hline
		2.5 & 1.2 & 26.3146 & 26.3146 & 4 & 1.3250 & 0.5000 & 0.8875  \\ \hline
		\end{tabular}
  \end{center}
\label{comp-hdp}
\end{table}
We compute the DSP and LSP for cubic loss function using same grid points for censoring time $\tau$ so that optimal sampling sampling plan in terms of $n$ and $\tau$ are same. We present the optimum sampling plans and the associated Bayes risks for different hyper parameters $a$ and $b$ in Table \ref{comp-hdp}.  In all the cases optimal DSP is as good as LSP in terms of Bayes risks. 
\par
 \noindent For non-polynomial loss function to obtain DSP, we consider the following values of coefficients:
$$ 
a_0 = 2,a_1 = 2,a_2 = 2, C_r = 30, C_s = 0.5, C_\tau = 0.5,  \xi^{*}=2, c^{*}=2.
$$ 
 and the results is given in Table \ref{NP-1}.
\begin{table}[H]
\caption{Minimum Bayes risk and corresponding optimal sampling plan for different values of $a$ and $b$ for DSP for non polynomial loss function }
	\vspace{1em}
	\begin{center}
		\begin{tabular}{|cc|ccccc|}\hline
		\multicolumn{1}{|c}{$a$} &
		\multicolumn{1}{c}{$b$} &
		\multicolumn{5}{|c|}{DSP} \\  
	  & & $r( n_0, \tau_0, \zeta_0, c_0)$ & $n_0$ & $\tau_0$ & $\xi_0$ & $c_0$  \\ \hline
    0.1 &    0.2 &    6.6966 &    2 &    0.6125 &    0.2250 &    1.6750\\ \hline
    1.0 &    0.2 &  26.1494  &   3 &    1.0875  &  0.3750 &   1.1500 \\ \hline
    1.5 &    0.8 &   19.4142 &   4 &    0.9000 &    0.3750 &    0.0750 \\ \hline
    2.5 &    0.8 &   27.5525 &    4 &    1.0625 &    0.3750 &   1.0875 \\ \hline
    3.0 &    0.8 &   29.6926 &   2 &    1.0750  &  0.3500 &    1.8250\\ \hline
    \end{tabular}
  \end{center}
\label{NP-1}
\end{table}
\subsubsection{\sc Comparison of DSP and LSP for  Type-I Hybrid censoring under quadratic loss function}
In this section, the numerical comparison between DSP and LSP is given for Type-I hybrid censoring. The values of coefficient $ a_0=2, a_1=2, a_2=2, C_s=0.5, r_s=0.3, C_\tau = 5.0$ and $ C_r=30$ are used for comparison. Table \ref{table-h1}-\ref{table-h2} represent the numerical results of comparison. The Bayes risk of LSP includes a complicated integrals which is computed by simulation techniques so Bayes risk of LSP is an approximation of exact Bayes risk of Bayesian sampling plan. 
\begin{table}[H]
\caption{Minimum Bayes risk and corresponding optimal DSP}
\centering
\vspace{1em}
		\begin{tabular}{|cc|c|cccccc|} \hline 
		\multicolumn{1}{|c}{$a$} &
		\multicolumn{1}{c|}{$b$} &
		\multicolumn{1}{c}{LSP \footnotemark} &
		\multicolumn{6}{|c|}{DSP}  \\ \cline{3-9}
	   & & $r(n_B,r_B,\tau_B,\delta_B)$  & $r( n_0, \tau_0, \xi_0,c_0)$ & $n_0$ & $r_0$ &  $\tau_0$ & $\zeta_0$ & $c_0$  \\ \hline
	   2.5 &   0.8 & 26.0319 & 26.0338   &  6   & 3  & 0.2000 & 0.2750 & 0.6600  \\ \hline
	   2.5  & 1.0  & 22.6430 & 22.6437 & 5 & 3 & 0.1875 &  0.2625 & 0.0725 \\ \hline
	   	3.0 & 0.8  & 29.7131 & 28.7890 & 4  & 2  & 0.2375  & 0.4250 & 0.0075 \\ \hline
	    \multicolumn{2}{|c}{$C_\tau$} & 	
	     \multicolumn{7}{|c|}{}  \\ \hline
		\multicolumn{2}{|c|}{0} & 24.6354 & 24.6736  &  4 &    4 &   0.8500  &  0.3000 & 0.3725 \\ \hline
	    \multicolumn{2}{|c|}{8} &  26.4662 & 26.4672 &   7 &    3 &    0.1625 &    0.2750 & 0.6600 \\ \hline
		\multicolumn{2}{|c|}{16} & 27.2513  & 27.2513 &    7 &    2 &    0.1000 &    0.2875 & 0.5775\\ \hline
	   	\end{tabular}
		 \label{table-h1}
\end{table}

\begin{table}[H]
\caption{Minimum Bayes risk and corresponding optimal DSP}
\centering
\vspace{1em}
		\begin{tabular}{|c|c|cccccc|} \hline 
		\multicolumn{1}{|c|}{} &
		\multicolumn{1}{c}{LSP \footnotemark[\value{footnote}]} &
		\multicolumn{6}{|c|}{DSP}  \\ \cline{2-8}
	   $C_s$ & $r(n_B,,r_B,\tau_B,\delta_B)$  & $r( n_0, \tau_0, \xi_0,c_0)$ & $n_0$ & $r_0$ &  $\tau_0$ & $\zeta_0$ & $c_0$  \\ \hline
	   0.5  & 26.0319 & 26.0338 & 6 & 3 & 0.2000 &  0.2750 & 0.6600  \\ \hline
	   0.6  & 26.5578 & 26.5626 & 5 & 3 & 0.2500 &  0.2750 & 0.6600 \\ \hline
	   	0.7  & 26.9106 & 26.9114 & 3  & 2  & 0.2750  & 0.3125 & 0.2400 \\ \hline
	    \multicolumn{1}{|c}{$C_r$} & 	
	     \multicolumn{7}{|c|}{}  \\ \hline
		\multicolumn{1}{|c|}{25} & 23.3583  & 23.3581 &  4 &  2  &   0.2375  &  0.3750 & 0.3350 \\ \hline
	    \multicolumn{1}{|c|}{30} & 26.0319 & 26.0338 &   6 &    3 &    0.2000 &    0.2750 & 0.6600 \\ \hline
		\multicolumn{1}{|c|}{40} &  30.0072 & 30.0071 &    7 &    4 &    0.1750 &    0.2375 & 0.1075\\ \hline
	   	\end{tabular}
		 \label{table-h2}
\end{table}
\footnotetext[1]{Bayes risk of LSP is obtain by simulation.}
Hence, for Hybrid Type-I censoring also the DSP is as good as LSP in terms of Bayes risk.

\section{\sc Conclusion}
In this paper we have worked on the improvement of the paper of \cite{Lam:1994} where he had used the MLE of the mean 
lifetime, which may not exist always for a Type-1 censored sample.  In fact \cite{LLH:2002} showed that the sampling plan proposed by \cite{Lam:1994} is neither optimal, 
nor Bayes.  The Bayesian sampling plan LSP proposed by  \cite{LLH:2002} provides a smaller Bayes risks 
than the sampling plan provided by \cite{Lam:1994}.   \cite{LLH:2002} implemented the LSP for quadratic loss function.  It is observed that in case of higher degree
polynomial loss function or for a more general non-polynomial loss function LSP may not be very easy to obtain.  In this
paper we have proposed a new decision theoretic sampling plan DSP based on an estimator which always exists and showed that it is as good as Bayesian sampling plan LSP in the sense that it minimizes the Bayes risk. It may be mentioned that although in this paper we have considered Type-I censored sample only but the method can be extended for other censoring cases also. More work is needed along that direction.

\section{\sc Acknowledgements} 
We express our sincere thanks to the referees and the editor for their useful suggestions which led to an improvement over an earlier version of this manuscript. 
\medskip

\section*{\sc Appendix}

To prove (\ref{brf}) first we show that

\be
 r(n, \tau, \xi, c)  =  nC_{s}+ \tau C_{\tau}+a_0 +a_1 \mu_1 +a_2 \mu_2 
+ \sum_{l=0}^{2}C_{l}\frac{b^a}{\Gamma (a)}\int_{0}^{\infty} \lambda^{a+l-1} 
                 e^{-\lambda b} P(\widehat{\theta}_{M+c} < \xi )d\lambda,   \label{bayes-risk}
\ee
\noindent {\sc Proof of} (\ref{bayes-risk}) 
\bea
 r(n, \tau, \xi, c) & = & E\big\{L(\delta(\textbf{x}),\lambda, n, \tau)\big\} \nonumber \\
            & = &   E_{\lambda}E_{X/\lambda}\big\{L(\delta(\textbf{x}),\lambda, n, \tau)\big\} \nonumber \\
            & = & 
            E_{\lambda}\big\{(nC_s + \tau C_{\tau} +a_0 +a_1 \lambda +a_2 \lambda^2)P(\hat{\theta}_{M+c} \geq \xi)  \nonumber \\
            &   & \quad +(nC_s + \tau C_{\tau} + C_r)P(\widehat{\theta}_{M+c} < \xi) \big\} \nonumber \\
            & = & nC_s + \tau C_{\tau} + a_0 + a_1 \mu_1 + a_2 \mu_2 \nonumber \\
            &  & \quad + E_{\lambda}\big\{(C_r - a_0 - a_1 \lambda - a_2 \lambda^2 )P(\widehat{\theta}_{M+c} < \xi)\big\}\nonumber \\
            & = & nC_s + \tau C_{\tau} + a_0 + a_1 \mu_1 + a_2 \mu_2 \nonumber \\
            &  & \quad + \int_0^{\infty} (C_r - a_0 - a_1 \lambda-a_2 \lambda ^2)
                \frac{b^a}{\Gamma (a) }\lambda^{a-1}e^{-\lambda b} P(\widehat{\theta}_{M+c} < \xi )\ d\lambda \nonumber \\
            & = & nC_{s}+ \tau C_{\tau}+a_0 +a_1 \mu_1 +a_2 \mu_2 
+ \sum_{l=0}^{2}C_{l}\frac{b^a}{\Gamma (a)}\int_{0}^{\infty} \lambda^{a+l-1} 
                 e^{-\lambda b} P(\widehat{\theta}_{M+c} < \xi )d\lambda.   \nonumber
\eea

\noindent {\sc Proof of} (\ref{brf})

Using (\ref{mix}) in (\ref{bayes-risk}) we get  
\bea
&  & \int_{0}^{\infty}\lambda^{a+l-1}e^{-\lambda b}P(\widehat{\theta}_{M+c}< \xi)\ d\lambda \nonumber \\
&  & = \int_{0}^{\infty}\lambda^{a+l-1}e^{-\lambda b}p S(\xi)\ d\lambda + \int_{0}^{\infty}\lambda^{a+l-1}e^{-\lambda b}(1-p)H(\xi)\ d\lambda \nonumber \\
&  & = \int_{0}^{\infty}\lambda^{a+l-1}e^{-\lambda b}e^{-n\lambda\tau} \ \textit{I}_{(\frac{n\tau}{c} < \xi)}\ d\lambda + \sum_{m=1}^{n}\sum_{j=0}^{m}\binom{n}{m}\binom{m}{j}(-1)^j \nonumber \\
&  &\hspace{3cm} \times \int_{0}^{\infty}\int_{0}^{\xi}\lambda^{a+l-1}e^{-\lambda \{b+\tau(n-m+j)\}} 
\pi \Big(x-\tau_{j;m,c}; m, \lambda(m+c)\Big)dx \ d\lambda     \nonumber\\
&  & = \frac{\Gamma{(a+l)}}{(b+n\tau)^{(a+l)}} \textit{I}_{(\frac{n\tau}{c} < \xi)} + \sum_{m=1}^{n}\sum_{j=0}^{m}\binom{n}{m}\binom{m}{j}(-1)^j \nonumber \\
&  &\hspace{3cm} \times \int_{0}^{\infty}\int_{0}^{\xi}\lambda^{a+l-1}e^{-\lambda \{b+\tau(n-m+j)\}} 
\pi \Big(x-\tau_{j;m,c}; m, \lambda(m+c)\Big)dx \ d\lambda     \nonumber \\
&  & =\frac{\Gamma{(a+l)}}{(b+n\tau)^{(a+l)}} \textit{I}_{(\frac{n\tau}{c} < \xi)} + \sum_{m=1}^{n}\sum_{j=0}^{m}\binom{n}{m}\binom{m}{j}(-1)^j\frac{(m+c)^m}{\Gamma (m)}\nonumber \\
&  & \hspace{4cm}\times \int_{0}^{\infty}\int_{\tau_{j;m,c}}^{\xi}\lambda^{a+l+m-1}e^{-\lambda\{b+(m+c)x \}}\big(x-\tau_{j;m,c}\big)^{m-1}dx \ d\lambda \nonumber \\
&  & =\frac{\Gamma{(a+l)}}{(b+n\tau)^{(a+l)}} \textit{I}_{(\frac{n\tau}{c} < \xi)} + \sum_{m=1}^{n}\sum_{j=0}^{m}\binom{n}{m}\binom{m}{j}(-1)^j\frac{(m+c)^m}{\Gamma (m)}\nonumber \\
&  & \hspace{7cm} \times \int_{\tau_{j;m,c}}^{\xi}\frac{\big(x-\tau_{j;m,c}\big)^{m-1}\Gamma{(a+l+m)}}{\{b+(m+c)x\}^{a+l+m}}dx \nonumber\\
&  & = \frac{\Gamma{(a+l)}}{(b+n\tau)^{(a+l)}} \textit{I}_{(\frac{n\tau}{c} < \xi)} + \sum_{m=1}^{n}\sum_{j=0}^{m}\binom{n}{m}\binom{m}{j}(-1)^j\frac{(m+c)^m}{\Gamma (m)} \nonumber \\  
&  & \hspace{5cm} \times \int_{0}^{\xi-\tau_{j;m,c}}\frac{y^{m-1}\Gamma{(a+l+m)}}{\{b+(m+c)\tau_{j;m,c}+(m+c)y\}^{a+l+m}}dy.    \label{step-2}    
\eea
Using $\ds C_{j,m}=b+(m+c) \tau_{j;m,c}$ in (\ref{step-2}), we can write
\begin{align}
\sum_{m=1}^{n}&\sum_{j=0}^{m}\binom{n}{m}\binom{m}{j}(-1)^j\frac{(m+c)^m}{\Gamma (m)}\frac{\Gamma{(a+l+m)}}{C_{j,m}^{a+l+m}}\int_{0}^{\xi-\tau_{j;m,c}}\frac{y^{m-1}}{\Big(1+\frac{(m+c)y}{C_{j,m}}\Big)^{a+l+m}}dy \nonumber\\
&=\sum_{m=1}^{n}\sum_{j=0}^{m}\binom{n}{m}\binom{m}{j}(-1)^j\frac{\Gamma{(a+l)}}{(C_{j,m})^{a+l}}\frac{\Gamma{(a+l+m)}}{\Gamma (m)\Gamma{(a+l)}}\int_{0}^{\frac{(m+c)(\xi- \tau_{j;m,c})}{C_{j,m}}}\frac{z^{m-1}}{(1+z)^{a+l+m}}dz.   \label{step-3}
\end{align}
Now taking a transformation $\ds z=u/(1-u)$, we have 
\begin{align*}
    \int_{0}^{C^*_{j,m,c}}\frac{z^{m-1}}{(1+z)^{a+l+m}}dz=\int_{0}^{S^*_{j,m,c}}u^{m-1}(1-u)^{a+l-1}du=B_{S^*_{j,m,c}}(m,a+l),
\end{align*} 
where $\ds C^{*}_{j,m,c}=\frac{(m+c)(\xi- \tau_{j;m,c})}{C_{j,m}}$, \ \ $\ds S^{*}_{j,m,c}=\frac{C^{*}_{j,m,c}}{1+C^{*}_{j,m,c}}$ , \hspace{0.5 cm} and 
\begin{equation}
 B_{x}(\alpha,\beta)=\int_{0}^{x}u^{\alpha -1}(1-u)^{\beta-1}du,  \ \ \ \ \ 0\leq x \leq 1, \nonumber
\end{equation}
is the incomplete beta function. Let us denote the cumulative distribution function of beta by
$$
I_x(\alpha,\beta)=B_x(\alpha,\beta)/B(\alpha,\beta).
$$ 
Then using (\ref{step-3}), (\ref{brf}) is finally obtained as
\begin{align}
r(n, \tau, \xi, c) = nC_{s}+ \tau C_{\tau}+a_0 + a_1 \mu_1 & +a_2 \mu_2  + \sum_{l=0}^{2}C_{l}\frac{b^a}{\Gamma (a)}\Bigg[\frac{\Gamma{(a+l)}}{(b+n\tau)^{(a+l)}} \textit{I}_{(\frac{n\tau}{c} < \xi)} 
 \nonumber \\            
 & + \sum_{m=1}^n\sum_{j=0}^m(-1)^{j}\binom{n}{m}\binom{m}{j} \frac{\Gamma{(a+l)}}{(C_{j,m})^{a+l}} 
I_{S^*_{j,m,c}}(m,a+l)\Bigg]. 
\end{align}

\noindent {\sc Proof of Result 3.1}

\begin{proof}
Note that the Bayes risk is given by 
\begin{align*}
\begin{split}
    r(n, \tau, \xi, c) & =E_{\lambda}\big\{(nC_s + \tau C_{\tau} + a_0 + a_1\lambda+ \ldots + a_k \lambda^k)P(\hat{\theta}_{M+c} \geq \xi) \\
               & + (nC_s + \tau C_{\tau} + C_r) P(\widehat{\theta}_{M+c} < \xi)\big\} \\
               & = nC_s + \tau C_{\tau} + E_{\lambda}\big\{(a_0 +a_1\lambda+ \ldots +a_k \lambda^k)P(\widehat{\theta}_{M+c} \geq \xi)\\
               & + C_r P(\widehat{\theta}_{M+c} < \xi)\big\}.
\end{split}
\end{align*}
Now we know that $a_0+a_1\lambda + \ldots +a_k\lambda^k \geq 0$ and $C_r$, the rejection cost, is non negative. Since 
$(n_0, \tau_0, \xi_0, c_0)$ is the optimal sampling plan so the corresponding  Bayes risk is 
\begin{align}
    r(n_0, \tau_0, \xi_0,c_0) \geq n_0 C_s + \tau_0 C_{\tau}.   \label{eq-14}
\end{align}
Now when $\xi=0$ we accept the batch without sampling and the corresponding Bayes risk is given by 
\begin{align*}
    r(0,0,0,0) = a_0+a_1\mu_1+ \ldots +a_k\mu_k.     \label{eq-15}
\end{align*}
When $ \xi=\infty$ we reject the batch without sampling and corresponding Bayes risk is given by 
\begin{align*}
    r(0,0,\infty,0) = C_r .
\end{align*}
Then the optimal Bayes risk is
\begin{align}
    r(n_0, \tau_0, \xi_0,c_0) \leq min \big\{ r(0,0,\infty,0),r(0,0,0,0),r(n, \tau, \xi_0(n,\tau), c_0(n,\tau))\big\}.
\end{align}
Hence from equations (\ref{eq-14}) and (\ref{eq-15}) we have
\begin{align*}
 n_0 C_s + \tau_0 C_{\tau} \leq min \big\{ r(0,0,\infty,0),r(0,0,0,0), r(n, \tau, \xi_0(n,\tau), c_0(n,\tau))\big\},
\end{align*}
from where it follows that 
\begin{align*}
n_0 \leq min \bigg\{ \frac{C_r}{C_s},\frac{a_0 + a_1\mu_1 + \ldots + a_k \mu_k}{C_s},\frac{r(n, \tau, \xi_0(n,\tau), c_0(n,\tau))}{C_s}\bigg\} \\
\tau_0 \leq min \bigg\{ \frac{C_r}{C_{\tau}},\frac{a_0 + a_1\mu_1 + \ldots + a_k \mu_k}{C_{\tau}},\frac{r(n, \tau, \xi_0(n,\tau), c_0(n,\tau))}{C_{\tau}}\bigg\} .
\end{align*}
\end{proof}

\nocite{*}
%\bibliographystyle{ieeetr}
%\bibliography{sample}
\section*{\sc References}
\addcontentsline{toc}{section}{Bibliography}
\small
\begingroup
\renewcommand{\section}[2]{}
% this must be set to use natbib (citep, citet) but requires BibTeX
\bibliographystyle{plainnat}

% number 99 determines how much citation can be included in file (maximum 99)

\end{document}